\documentclass[11pt]{amsart}

\usepackage{multirow}
\usepackage{graphicx}
\usepackage[top=26mm, bottom=26mm, left=26mm, right=26mm]{geometry}
\usepackage[dvips]{color}
\usepackage{bm}
\usepackage{marvosym}
\theoremstyle{definition}
\newtheorem{remark}{Remark}
\numberwithin{remark}{section}

\numberwithin{definition}{section}
\newcommand{\vct}[1]{\bm{\mathsf{#1}}}

\newcommand{\mtx}[1]{\bm{\mathsf{#1}}}

\newcommand{\vtwo}[2]{\left[\begin{array}{cc} #1 \\ #2 \end{array}\right]}

\newcommand{\pgnotate}[1]{}

\newcommand{\lsp}{\vspace{3mm}}

\newcommand{\bi}{\begin{itemize}}

\newcommand{\ei}{\end{itemize}}

\newcommand{\ben}{\begin{enumerate}}

\newcommand{\een}{\end{enumerate}}

\newcommand{\be}{\begin{equation}}

\newcommand{\ee}{\end{equation}}

\newcommand{\bea}{\begin{eqnarray}}

\newcommand{\eea}{\end{eqnarray}}

\newcommand{\ba}{\begin{align}}

\newcommand{\ea}{\end{align}}

\newcommand{\bse}{\begin{subequations}}

\newcommand{\ese}{\end{subequations}}

\newcommand{\bc}{\begin{center}}

\newcommand{\ec}{\end{center}}

\newcommand{\bfi}{\begin{figure}}

\newcommand{\efi}{\end{figure}}

\newcommand{\bmp}[1]{\begin{minipage}{#1}}

\newcommand{\emp}{\end{minipage}}

\newcommand{\bp}{\begin{proof}}

\newcommand{\ep}{\end{proof}}

\begin{document}

\begin{center}An integral equation technique for scattering problems with mixed boundary conditions
\textbf{\large }

\lsp

 Adrianna Gillman\\
Department of Computational and Applied Mathematics\\
Rice University

\end{center}
\lsp

\noindent
\textbf{Abstract:}  This paper presents an integral formulation for Helmholtz 
 problems with mixed boundary conditions.  Unlike most integral equation techniques 
 for mixed boundary value problems, the proposed method uses a global boundary charge 
 density.  As a result, Calder\'on identities can be utilized to avoid the use of 
 hypersingular integral operators.  More importantly, the formulation avoids spurious 
 resonances.  Numerical results illustrate the performance of the proposed solution technique.
%  A comparison of the performance of the proposed method and the integral formulation of \cite{2009_Helsing_mixedBVP} 
%  for Laplace mixed boundary value problems is also presented.\end{abstract}

\lsp

\section{Introduction}

This manuscript describes an integral equation technique for solving Helmholtz 
mixed boundary value problems for a scattering body $\Omega$ with boundary $\Gamma = \partial \Omega$.  
Specifically, we consider the following boundary value problem

\begin{equation}
\label{eq:mixedBVP}
% \tag{DBVP}
\left\{
\begin{aligned}
\Delta u + \omega^2 u =&\ 0\qquad&&\mbox{ in } \Omega^c,\\
                    u =&\ g\qquad&&\mbox{ on } \Gamma_D,\\
                    \frac{\partial u}{\partial \nu} =&\ f\qquad&&\mbox{ on } \Gamma_N,\\
%                     \frac{\partial u}{\partial \nu}+i\xi u =&\ h\qquad&&\mbox{ on } \Gamma_R,\\
%                     \frac{\partial u} {\partial \nu} +i\eta\u  = h& \qquad && \mbox{ on } \Gamma_R,\\
\frac{\partial u}{\partial r}  - i\omega  u  =&\ O  \left( \frac{1}{r} \right)\qquad&&\mbox{ as } r\to\infty,
\end{aligned}\right.
\end{equation}
where $\Gamma_D$ denotes the portion of $\Gamma$ with Dirichlet boundary condition, $\Gamma_N$ 
denotes the portion of $\Gamma$ with Neumann boundary condition
% , $\Gamma_R$ denotes the portion 
% of the boundary with (impedance) Robin boundary condition 
and $\Gamma= \Gamma_D\cup\Gamma_N$. 
The last equation in (\ref{eq:mixedBVP}) is the outgoing Sommerfeld radiation condition which specifies 
a decay condition on the solution.

\subsection{Prior work}  
Several integral equation based solution techniques have 
previously been used for solving mixed boundary value problems.  

In \cite{2009_Helsing_mixedBVP,2012_Ojala_mixedBVP}, Laplace mixed boundary value problems are considered. 
Both papers represent the solution as a linear combination of a double layer integral operator 
defined on $\Gamma_D$ with a single layer integral operator defined on 
$\Gamma_N$.  The integral equation that results from enforcing boundary 
conditions involves evaluating a hypersingular integral operator which has 
a norm that can grow without bound even when considered on neighboring segments of $\Gamma$.  
By utilizing a near evaluation technique designed for the hypersingular operator and local
compression, there does not appear to be a loss of 
accuracy associated with this solution technique.  

The reformulation of Helmholtz boundary value 
problems as boundary integral equations needs to be done with 
some care to avoid introducing spurious 
resonances.  A {\it{spurious resonance}} 
is a wavenumber $\omega$ for which a na\"ive integral 
equation formulation would result in a non-trivial null space even when the 
original boundary value problem is well-posed.

\cite{Kress77} considers (\ref{eq:mixedBVP}) where $\Omega$ consists of a 
collection of scattering bodies and each body has a single type 
of boundary condition.  To avoid spurious resonances, \cite{Kress77} uses a linear combination of 
combined field representations (one for $\Gamma_D$ and another for $\Gamma_N$).
Roughly speaking, a block preconditioner 
involving two inverses is utilized to achieve a second kind integral equation 
for (\ref{eq:mixedBVP}).  Unfortunately, this technique does not extend 
directly to single body mixed boundary condition problems.  \cite{Kleefeld} presents a 
similar approach for solving (\ref{eq:mixedBVP}) on a single scattering body which
utilizes a preconditioner that is equivalent to squaring the inverse of a first kind integral equation.  

At this point, there does not appear to exist a uniquely solvable pure second kind integral equation for (\ref{eq:mixedBVP}).

\subsection{Outline}
This paper presents an integral formulation for solving (\ref{eq:mixedBVP}) which
uses a single regularized combined field integral representation for the solution.
Utilizing a Calder\'on identity leads to a simple block first kind integral equation system.  

The manuscript begins by reviewing robust 
integral formulations for Dirichlet and Neumann boundary value problems in section  \ref{sec:formulations}.
Then the proposed integral formulation is presented.  Next, section \ref{sec:discrete}
presents a technique for discretizing the resulting integral equation system.  
Finally, section \ref{sec:numerics} illustrates the performance of 
the proposed solution technique. 

\begin{remark}
 While this paper focuses on Helmholtz problems with mixed Dirichlet and Neumann boundary data,
 the solution technique can easily be extended to case of additional Robin boundary conditions.  Also,
 Laplace boundary value problems can be handled by setting $\omega =0$ and using the appropriate 
 Green's function.
\end{remark}

%%%%%%%%%%%%%%%%%%%%%%%%%%%%%%%%%%%%%%%%%%%%%%%%%%%%%%%%%%%%%%%%%%%%%%%%%%%%%%%
\section{Integral equation techniques}
\label{sec:formulations}
The reformulation of (\ref{eq:mixedBVP}) as a block integral equation system 
involves the classical single, hypersingular and double layer kernels

\begin{equation}
\begin{split}
S_\omega f(x)   &=   \int_{\Gamma} G_\omega (x,y) f(y) ds(y),                        \\
T_\omega f(x) &= \eta_x \cdot \nabla_x\int_\Gamma \eta_y \cdot \nabla_y G_\omega (x,y) f(y) ds(y), \\
D_\omega f(x)   &=   \int_{\Gamma} \eta_y \cdot \nabla_y G_\omega (x,y) f(y) ds(y)\ \ \mbox{and} \\
D_\omega ^*f(x) &=   \int_{\Gamma} \eta_x \cdot \nabla_x G_\omega (x,y) f(y) ds(y)
\end{split}
\label{formulations:operators}
\end{equation}
where $G_\omega (x,y)$ is the free space Green's function
$$G_\omega (x,y) = \frac{\rm{i}}{4}H_0(\omega |x-y|)$$
of the two dimensional Helmholtz problem with wavenumber $\omega$  
and $\eta_p$ denotes the outward facing normal vector at the point $p\in \Gamma$.

This section begins by presenting robust integral formulations for Dirichlet and 
Neumann boundary value problems.  Then we propose a block integral equation system to solve (\ref{eq:mixedBVP}).

\subsection{The Dirichlet boundary value problem}

Consider the Dirichlet boundary value problem
\begin{equation}
\label{eq:Dir}
% \tag{DBVP}
\left\{
\begin{aligned}
\Delta u + \omega^2 u =&\ 0\qquad&&\mbox{ in } \Omega^c,\\
                    u =&\ g\qquad&&\mbox{ on } \Gamma,\\
%                     \frac{\partial u}{\partial \nu} =&\ f\qquad&&\mbox{ on } \Gamma_N,\\
%                     \frac{\partial u}{\partial \nu}+i\xi u =&\ h\qquad&&\mbox{ on } \Gamma_R,\\
%                     \frac{\partial u} {\partial \nu} +i\eta\u  = h& \qquad && \mbox{ on } \Gamma_R,\\
\frac{\partial u}{\partial r}  - i\omega  u  =&\ O  \left( \frac{1}{r} \right)\qquad&&\mbox{ as } r\to\infty,
\end{aligned}\right.
\end{equation}
where $\omega$ denotes the constant wavenumber.

The solution can be represented as a combined field
\begin{equation}
\label{eq:combine}
u(x) = D_\omega  \sigma(x) - i|\omega| S_\omega  \sigma(x),
\end{equation}
where $\sigma(x)$ is an unknown boundary charge density.  

It is well-known \cite{Nedelec,Colton-Kress2} that a uniquely solvable boundary integral
 equation 
 \begin{equation}
\frac{1}{2} \sigma(x)  + D_\omega \sigma(x) - i |\omega|  S_\omega  \sigma(x) = g(x)
\label{formulations:combined}
\end{equation}
for $\sigma(x)$ results from enforcing the boundary condition.  The integral 
equation (\ref{formulations:combined}) 
is second kind on smooth $\Gamma$ and is sometimes referred to as 
the {\it{Combined Field Integral Equation}}.

\subsection{The Neumann boundary value problem}
\label{sec:Neu}

Consider the Neumann boundary value problem
\begin{equation}
\label{eq:Neu}
% \tag{DBVP}
\left\{
\begin{aligned}
\Delta u + \omega^2 u =&\ 0\qquad&&\mbox{ in } \Omega^c,\\
%                     u =&\ g\qquad&&\mbox{ on } \Gamma,\\
                    \frac{\partial u}{\partial \nu} =&\ f\qquad&&\mbox{ on } \Gamma,\\
%                     \frac{\partial u}{\partial \nu}+i\xi u =&\ h\qquad&&\mbox{ on } \Gamma_R,\\
%                     \frac{\partial u} {\partial \nu} +i\eta\u  = h& \qquad && \mbox{ on } \Gamma_R,\\
\frac{\partial u}{\partial r}  - i\omega  u  =&\ O  \left( \frac{1}{r} \right)\qquad&&\mbox{ as } r\to\infty,
\end{aligned}\right.
\end{equation}
where $\omega$ denotes the constant wavenumber, and $\nu$ denotes the outward facing normal vector for $x\in\Gamma$.

Using the combined field representation 
(\ref{eq:combine}) for the solution of (\ref{eq:Neu}) results in the integral equation 

\begin{equation}
 -i|\omega|\left(-\frac{1}{2}\sigma(x)+ D^*_\omega\sigma(x)\right) +T_\omega \sigma(x) = f(x)
 \label{eq:Neucombined}
\end{equation}

which is not a second kind integral equation for smooth $\Gamma$ since the 
hypersingular integral operator $T_\omega$ is not compact.  The $T_\omega $ operator is 
troublesome for two reasons.  First, $T_\omega$ is an unbounded operator from $L^2(\Gamma)$ to 
$L^2(\Gamma)$.  As a direct consequence, the linear system 
resulting from discretization of the integral equation (\ref{eq:Neucombined}) is ill-conditioned.
Second, standard singular quadrature such as \cite{2001_rokhlin_kolm} are not 
sufficient to discretize the operator.  While it is possible to overcome 
both these difficulties by developing analytical and quadrature techniques which view $T_\omega$ as an
operator between two Sobolev spaces, there are highly accessible alternatives.  

 The most common approach to avoid these problems
is to utilized a so-called {\it{regularized}} combined field representation 
\begin{equation}
u(x) = D_\omega S_\omega \sigma(x) - i |\omega|  S_\omega  \sigma(x)
 \label{eq:reg_comibine}
\end{equation}
where $\sigma(x)$ still represents an unknown boundary charge distribution.  
Enforcing the Neumann boundary condition results in a regularized boundary integral equation
\begin{equation}
 -i|\omega|\left(-\frac{1}{2}\sigma(x)+ D^*_\omega\sigma(x)\right) +T_\omega S_\omega \sigma(x) = f(x).
 \label{eq:NeuREGcombined}
\end{equation}
This integral equation is called regularized because upon utilizing Calder\'on identities \cite{Nedelec}
we can rewrite (\ref{eq:NeuREGcombined}) as
\begin{equation}
\left(\frac{1}{4}+i|\omega|\frac{1}{2}\right)\sigma(x)-i|\omega|D^*_\omega\sigma(x) +(D^*_\omega)^2 \sigma(x) = f(x).
\label{eq:CALREGcombined}
\end{equation}
which does not involve hypersingular integral operators.

% An alternative technique is to precondition (\ref{eq:Neucombined}) by $S_\omega$ and utilize a different Calder\'on identity.
% The result is a different integral equation 
% \begin{equation}
%  \frac{1}{4}\sigma(x)-D_\omega^2\sigma(x)-i|\omega|S_\omega\left(-\frac{1}{2}\sigma(x)+D^*_\omega\sigma(x) \right)= S_\omega f(x)
%  \label{eq:precond}
% \end{equation}
% which also does not involve hypersingular integral operators.

On smooth geometries, equation (\ref{eq:CALREGcombined}) is 
a second kind integral equation.

\subsection{The mixed boundary value problem}
For the mixed boundary value problem (\ref{eq:mixedBVP}), we choose to 
use the regularized combined field representation (\ref{eq:reg_comibine}) 
for the solution.  With this choice, we are able to utilize to two very appealing properties: 
(i) it avoids spurious resonances and (ii) special quadrature for hypersingular integral
operators is not needed.  The compromise is that the resulting integral equation has 
a first kind block. 

From section \ref{sec:Neu}, we know the block row integral equation for $x\in\Gamma_N$ 
is given by equation (\ref{eq:NeuREGcombined}).  

Applying (\ref{eq:reg_comibine}) to $\Gamma_D$ results in a first kind integral equation given by 
\begin{equation}\label{eq:DirREGcombined}\frac{1}{2}S_\omega\sigma(x)+D_\omega S_\omega\sigma(x)-i|\omega|S_\omega\sigma(x) = f(x).\end{equation}

Rewriting the integral equation in block form, we get

\begin{equation}
\vtwo{S^D_\omega +(D_\omega S_\omega)^D-i|\omega|S_\omega^D}{(\frac{1}{4}+i|\omega|\frac{1}{2})I^N -i|\omega|D^{*,N}_\omega +(D^*_\omega)^{2,N}}
% \left(\mtwo{S^{DD}_{\omega}}{\vct{0}}{\vct{0}}{(\frac{1}{4}+i|\omega|\frac{1}{2})\mtx{I}}
%  +
% \mtwo{(D_\omega S_{\omega})^{DD}-i|\omega|S_\omega^{DD}}
% {(D_\omega S_{\omega})^{DN}-i|\omega|S_\omega^{DN}}
% {-i|\omega|D^{*,ND}_\omega +(D^*_\omega)^{2,ND}}{-i|\omega|D^{*,NN}_\omega +(D^*_\omega)^{2,NN}}\right)
% \vtwo{\sigma_D}{\sigma_N} 
\sigma
= \vtwo{g}{f}
 \label{eq:sys}
\end{equation}
where  
$$S_{\omega}^{D}\sigma(x) = \int_{\Gamma}G_\omega(x,y)\sigma(y)ds(y) \ {\rm for} \ x\in\Gamma_D,$$ 
$$D^{*,N}_{\omega}\sigma(x) = \int_{\Gamma}G_\omega(x,y)\sigma(y)ds(y) \ {\rm for}  \ x\in\Gamma_N,$$
$$(D_\omega S_\omega)^{D}\sigma(x) = \int_\Gamma \eta_y \cdot \nabla_y G_\omega (x,y) \left(\int_{\Gamma} G_\omega(y,w)\sigma(w)ds(w)\right)ds(y) \ x\in\Gamma_D$$
etc.

% We write (\ref{eq:sys}) in condensed form as
% $$\mtx{A}\vct{\sigma}=\vct{b}.$$
Let $\mtx{A}\vct{\sigma} =\vct{b}$ denote the condensed form of (\ref{eq:sys}).

% Hence we choose to represent the solution with the combined field (\ref{eq:combine}).
% Then we have to solve a system of equations corresponding to the restriction 
% of (\ref{formulations:combined}) for $x\in\Gamma_D$ and (\ref{eq:precond}) for $x\in\Gamma_N$.  
% 
% 
% 
% To clarify, for the mixed boundary value problem (\ref{eq:mixedBVP}), we choose 
% to represent the solution as
% \begin{equation*}
% % \label{eq:combine}
% u(x) = D_\omega  \sigma(x) - i|\omega| S_\omega  \sigma(x).
% \end{equation*}
% Enforcing the boundary conditions results in block integral equation system
% \begin{equation}\frac{1}{2} \sigma(x)  + D_\omega \sigma(x) - i |\omega|  S_\omega  \sigma(x) = g(x) \qquad {\rm for} \ x\in \Gamma_D\label{eq:dir}\end{equation}
% and 
% \begin{equation} \frac{1}{4}\sigma(x)-D_\omega^2\sigma(x)-i|\omega|S_\omega\left(-\frac{1}{2}\sigma(x)+D^*_\omega\sigma(x) \right)= S_\omega f(x) 
% \qquad {\rm for} \ x \in \Gamma_N.\label{eq:neu}\end{equation}
% 

%%%%%%%%%%%%%%%%%%%%%%%%%%%%%%%%%%%%%%%%%%%%%%%%%%%%%%%%%%%%%%%%%%%%%%%%%%%%%%%
\section{Discretization}
\label{sec:discrete}

\begin{figure}[h]
 \centering
\setlength{\unitlength}{1mm}
 \begin{picture}(180,60)
  \put(-05,05){ \includegraphics[width=75mm]{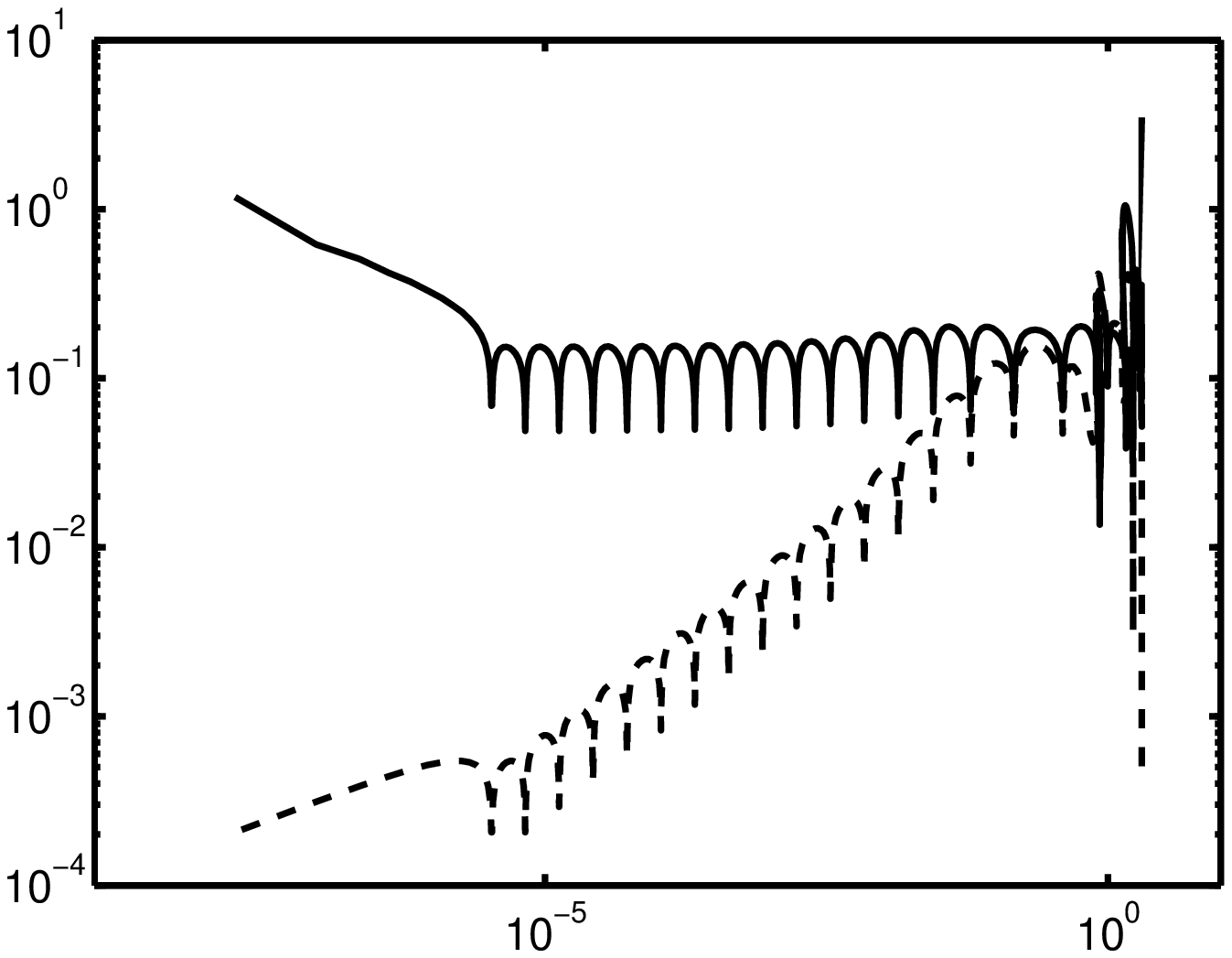}}
  \put(-05,30){\rotatebox{90}{$|\sigma|$}}
  \put(35,05){$r$}
  \put(33,00){(a)}
  \put(75,05){ \includegraphics[width=75mm]{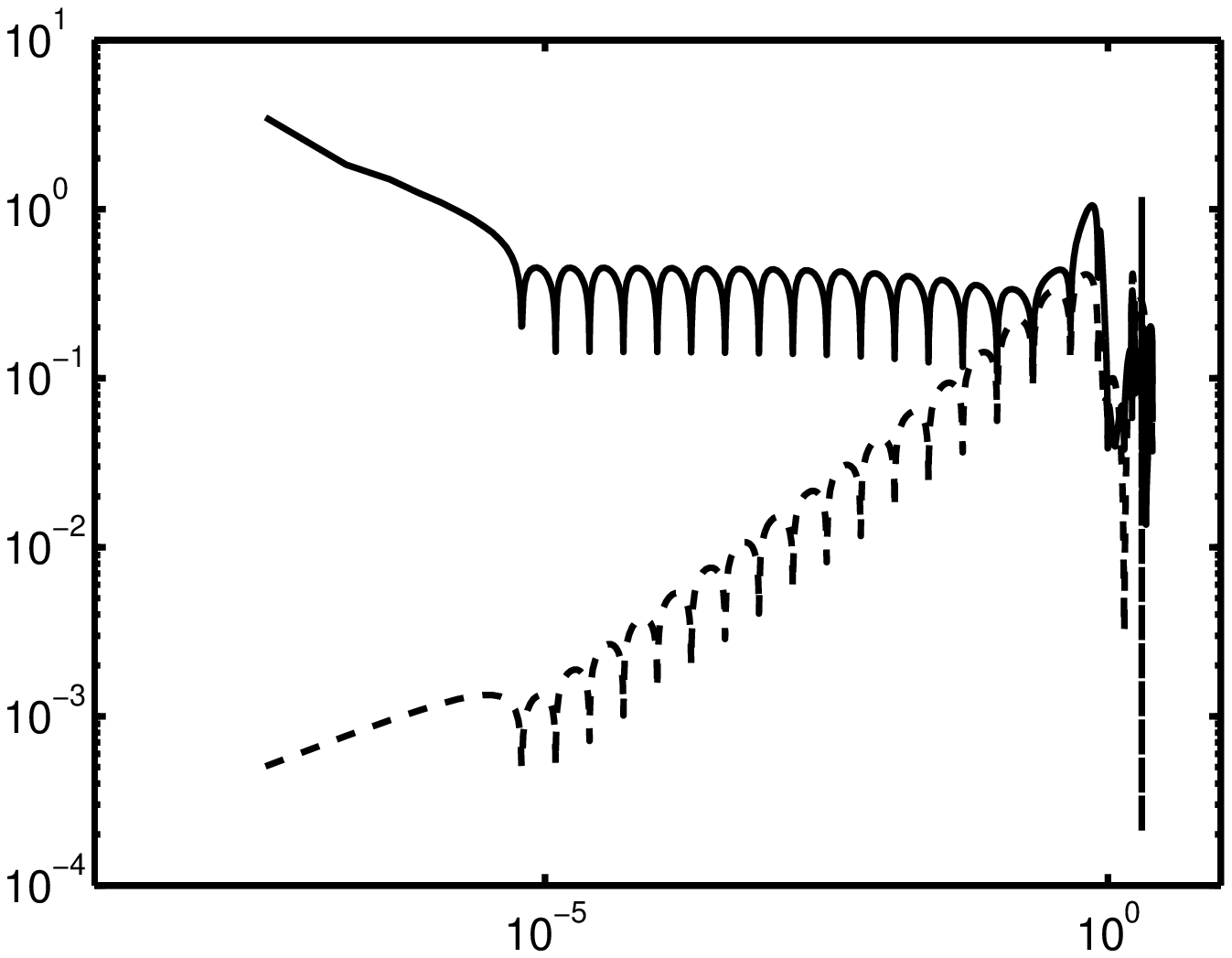}}
  \put(75,30){\rotatebox{90}{$|\sigma|$}}
  \put(115,05){$r$}
  \put(113,00){(b)}
 \end{picture}
\caption{\label{fig:densities} Illustration of $|\sigma|$ on $\Gamma_D$ (solid line) and 
$\Gamma_N$ (dashed line) vs the distance $r$ from the Dirichlet-Neumann 
Junction on the left-hand side (a) and right-hand side (b) of the star geometry (see Figure \ref{fig:geom})(a)
when $\omega = 1$ and the solution is unknown (see section \ref{sec:helm} for problem specifications).  
Four levels of dyadic refinement are utilized making the closest points to the junctions $r = 10^{-8}$ away.}
\end{figure}

Some care needs to be taken when discretizing the block integral equation (\ref{eq:sys})
since the boundary charge distribution $\sigma$ is likely to not be smooth 
at the Dirichlet-Neumann junctions.   Figure \ref{fig:densities} illustrates the behavior
of the boundary charge distribution $\sigma$ when $\omega =1$ and half a smooth star geometry 
(see Figure \ref{fig:geom}(a)) has zero Neumann boundary 
condition while the remainder has a Dirichlet boundary condition generated by a solution a free space 
Helmholtz problem.
We follow the approach of 
\cite{2012_bremer_direct_BIE_corners,2011_helsing_corner_BIE,2008_helsing_corner_BIE,2011_helsing_squares,2014_corner_note}
which describe techniques for dealing with the singularities that occur 
when solving scattering problems on Lipschitz geometries.  Thus a standard panel based quadrature 
rule is used for the Nystr\"om discretization of (\ref{eq:sys}).  
Near the Dirichlet-Neumann junctions the mesh is refined until the contribution from the 
panels closest to the junctions is small.

\begin{figure}[h]
 \centering
 \begin{tabular}{ccccc}
  \includegraphics[width=50mm]{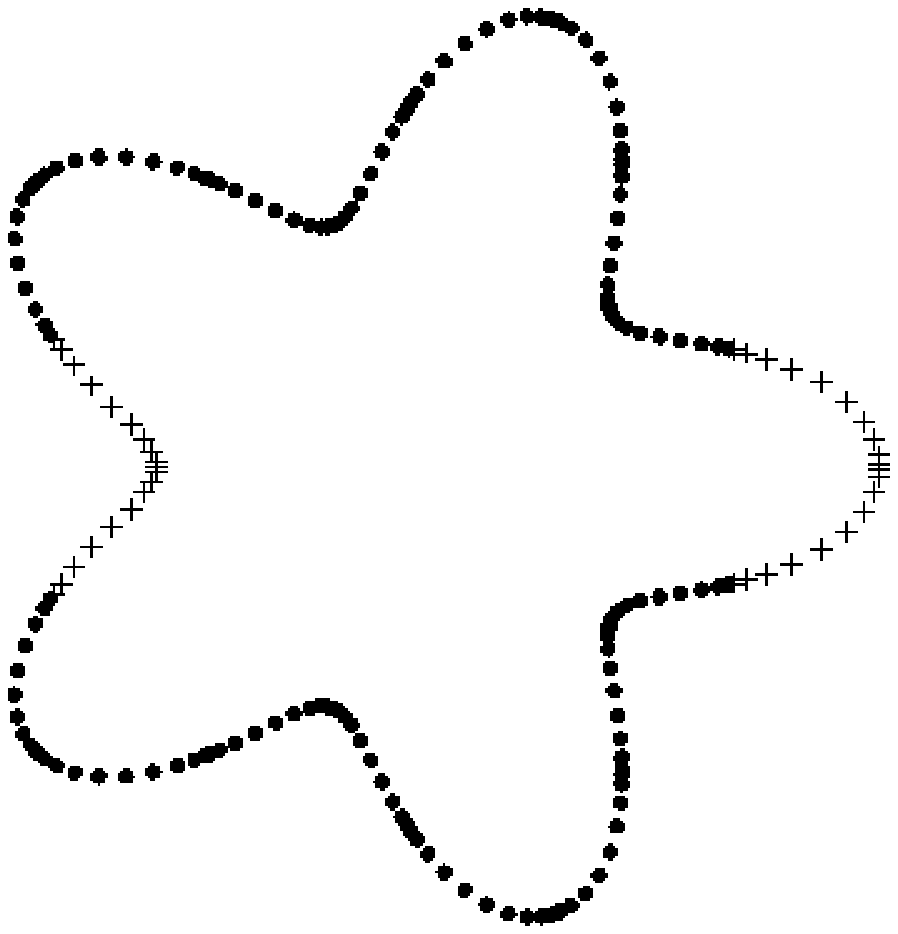} &\mbox{\hspace{3mm}} &\includegraphics[width=50mm]{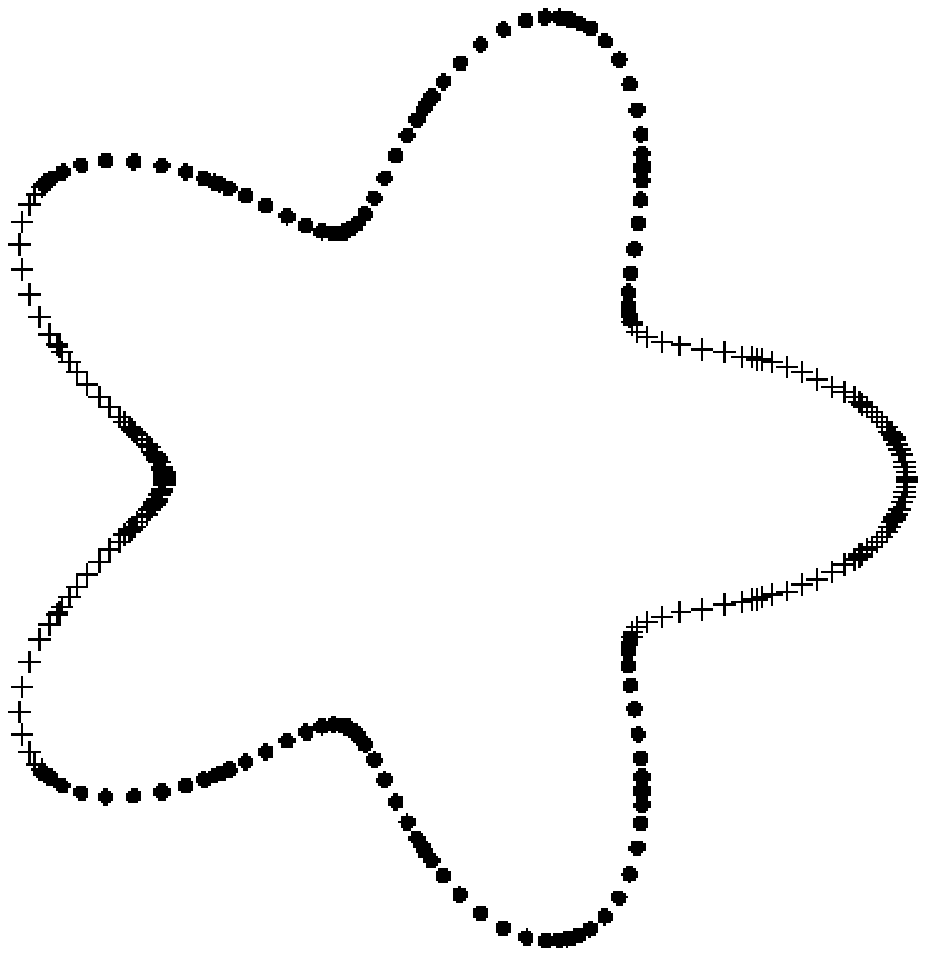} &\mbox{\hspace{3mm}}& \includegraphics[height=50mm]{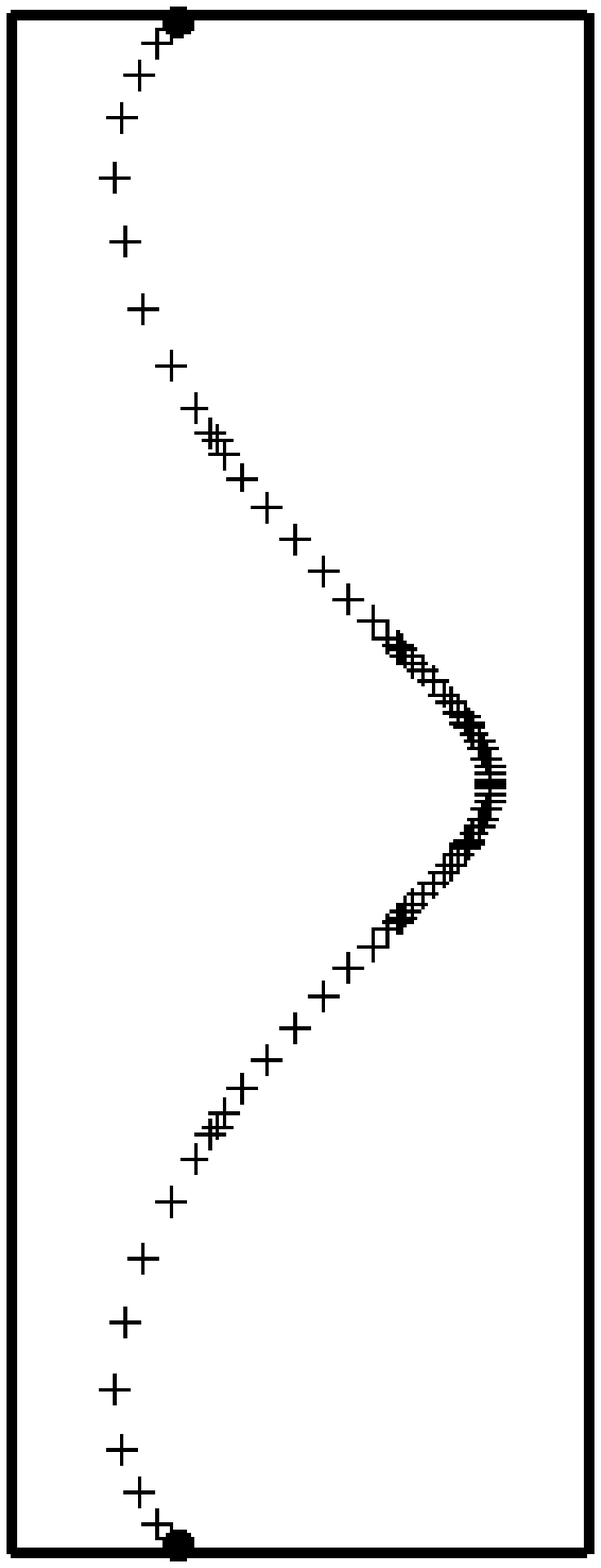}\\
  (a) &\mbox{\hspace{3mm}}& (b) &\mbox{\hspace{3mm}}& (c)
  \end{tabular}
\caption{\label{fig:discretization} Illustration of the discretization for 
a smooth star geometry 
with $\Gamma_N$ corresponding to the upper half of $\Gamma$ and $\Gamma_D$ corresponding to the
lower half of $\Gamma$ where $\bullet$ are nodes in $\Gamma_s$ and $+$ are nodes in $\Gamma_r$.  
(a) A 10 point composite Gaussian grid with no refinement. (b) The three level locally refined grid. (c) A closer 
view of the locally refined grid.   } 
\end{figure}

For simplicity of presentation, consider a geometry where $\Gamma_D$ is 
comprised solely of one section of $\Gamma$ and (by default) 
$\Gamma_N$ is also one portion of $\Gamma$ as in Figure \ref{fig:geom}.  First the geometry is partitioned into 
the six pieces such that  $\Gamma_1\cup\Gamma_2\cup\Gamma_3 = \Gamma_D$ and 
$\Gamma_4\cup\Gamma_5\cup\Gamma_6 = \Gamma_N$.  $\Gamma_2$ and $\Gamma_5$ are 
the parts of $\Gamma$ not directly touching the boundary condition junctions, while the pairs
$\{\Gamma_1,\Gamma_6\}$ and $\{\Gamma_3,\Gamma_4\}$ join at boundary condition junctions.
Hence $\Gamma_s = \Gamma_2\cup\Gamma_5$ is the portion of $\Gamma$ where the 
solution $\sigma$ is smooth and 
$\Gamma_r = \Gamma_1\cup\Gamma_3\cup\Gamma_4\cup\Gamma_6$ is the portion of $\Gamma$ 
where local refinement is likely needed. 
Figure \ref{fig:discretization}(a) illustrates the partitioning of 
$\Gamma$ for the smooth star geometry where the 
upper half of $\Gamma$ has Neumann boundary conditions while the lower half of $\Gamma$ 
has Dirichlet boundary conditions.  Since $\sigma$ is smooth on $\Gamma_s$, this region can be 
discretized coarsely. For the portions of the boundary near the boundary 
condition junctions $\Gamma_r$, the panels nearest 
to the junctions are recursively cut in half until the contribution from the panel nearest the junction 
is so small that it can be ignored. For many problems on smooth geometries, it is sufficient to 
stop $10^{-8}$ away from the junction point to obtain eight digits of accuracy.  
Referencing the third column of table \ref{tab:star_helm}, we note that using four levels of refinement 
corresponding to being approximately $10^{-8}$ away from the junctions does result in 
eight digits of accuracy.

Figure \ref{fig:discretization}(b) illustrates the 
discretization with three levels of refinement into the junctions. Figure \ref{fig:discretization}(c)
is a zoomed in view of the left-hand side of the same figure. Note that this refinement 
procedure introduces a superfluous amount of points near the junctions.  
The extra points can be eliminated via compression techniques presented in 
\cite{2012_bremer_direct_BIE_corners,2011_helsing_corner_BIE,2008_helsing_corner_BIE,2011_helsing_squares,2014_corner_note}.
% \agnotate{verify that all of these papers reference the correspond to compression.}
Since the focus of this paper is on the performance of the integral formulation, no compression techniques are 
utilized in the numerical experiments in section \ref{sec:numerics}.

% For problems where the number of required discretization points is small, equations (\ref{eq:NeuREGcombined}) and (\ref{eq:DirREGcombined})
% can be made for the full geometry then the discretized form of (\ref{eq:sys}) is found by looking at the restriction of those 
% systems to the appropriate portions of the geometry.    For problems requiring a large number of discretization points, more care 
% is used to 

%%%%%%%%%%%%%%%%%%%%%%%%%%%%%%%%%%%%%%%%%%%%%%%%%%%%%%%%%%%%%%%%%%%%%%%%%%%%%%%%%
\section{Numerical examples}
\label{sec:numerics}

This section reports on the performance of the purposed solution technique for solving boundary value problems 
with mixed boundary data.  Section \ref{sec:helm} considers the Helmholtz problem with mixed boundary conditions for 
three geometries, a smooth star, a pacman and a tear geometry.  These geometries represent the three most common 
boundary phenomena:  smooth boundary, a re-entrant corner, and a corner with the possibility of 
confined oscillations.  Performance of the solution technique for the boundary conditions corresponding to 
a smooth known solution and an unknown solution are presented.  Section \ref{sec:Lap} considers a Laplace problem 
with mixed boundary conditions on a smooth geometry.  A comparison of the performance of the proposed method 
versus the integral formulation proposed in \cite{2009_Helsing_mixedBVP} is reported.  

All integral equations are discretized using a Nystr\"om
technique based on a 16-point composite Gaussian quadrature \cite{2001_rokhlin_kolm}.   Panels are 
placed on the geometries via the method described in section \ref{sec:discrete}.  The 
number of discretization points is $N = 16(N_{\rm pan}+4(2^l))$ where $N_{\rm pan}$ 
is the number of panels on $\Gamma_s$ and $2^l$ is the number of refinement panels 
utilized.

The experiments were run on a Lenovo laptop computer with 16GB of RAM and a 2.4GHz Intel i7-4700M procesor
in Matlab.

\subsection{Helmholtz problems}
\label{sec:helm}

\begin{figure}[h]
\centering
\begin{tabular}{ccc}
\includegraphics[height=4cm]{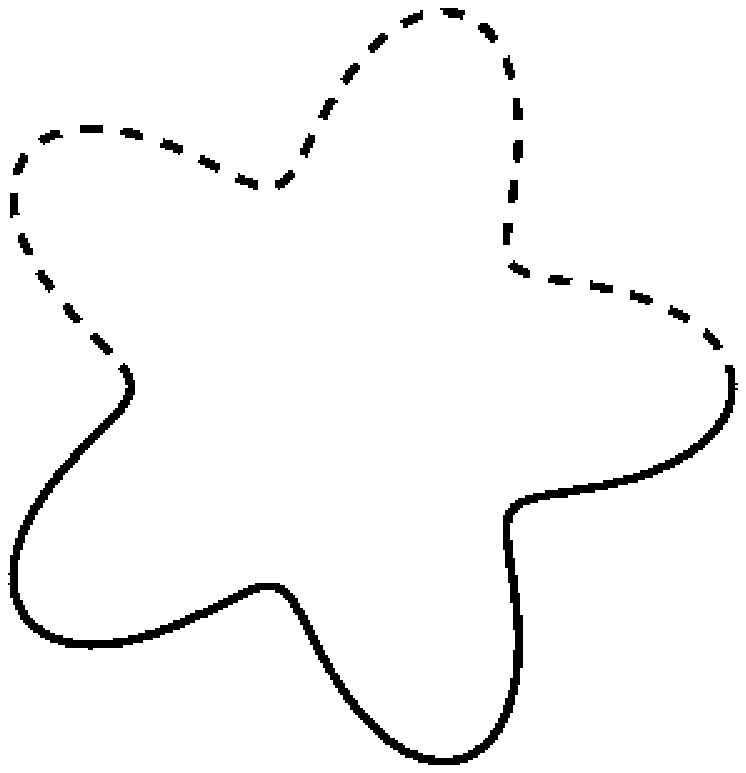}  &
\includegraphics[height=3.8cm]{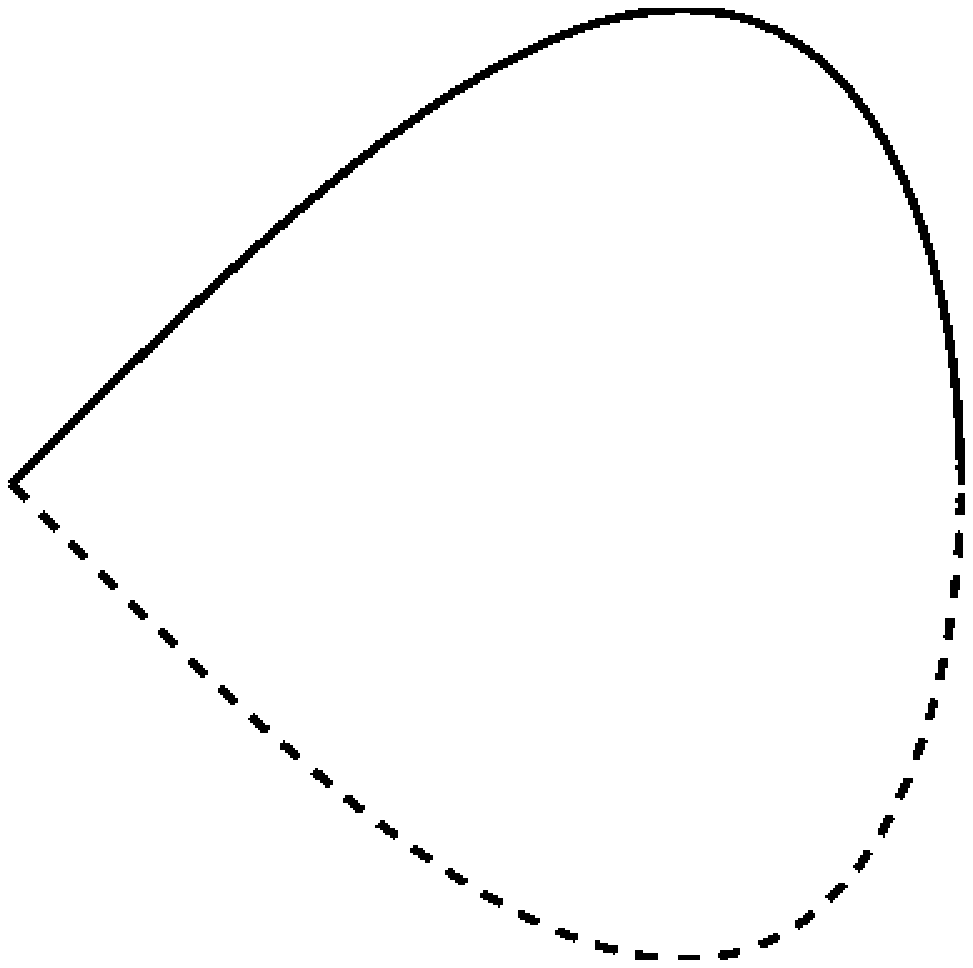} &
\includegraphics[height=4cm]{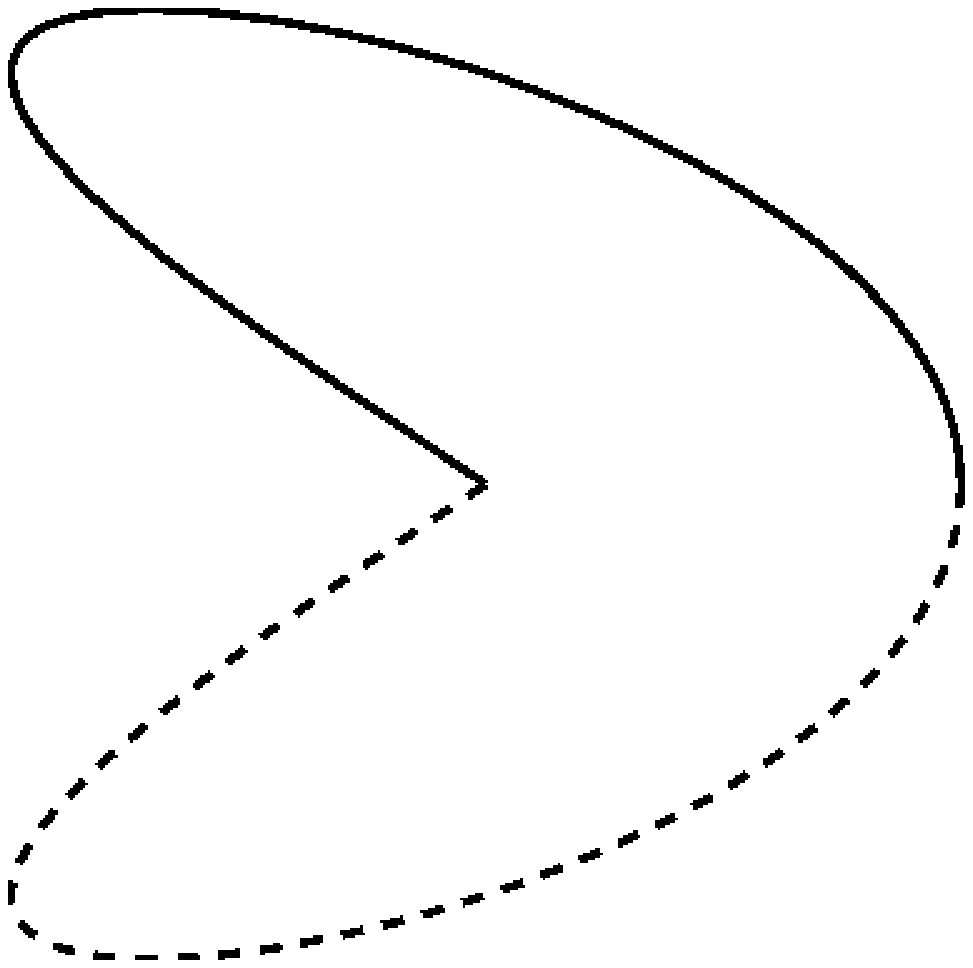} \\
(a) & (b) & (c)
\end{tabular}
\caption{\label{fig:geom} The contours $\Gamma$ used in the numerical
experiments in Section \ref{sec:numerics}. (a) Smooth star. (b) Tear. (c) Pacman. The portion 
of the boundary with the dashed line has Neumann boundary conditions while the portion of the 
boundary with solid line has Dirichlet boundary conditions.}
\end{figure}

This section considers the mixed boundary value problem (\ref{eq:mixedBVP}) with three different wave numbers 
($\omega = 1$, $10$, and $100$) on three geometries: 

\begin{list}{}{}
 \item[smooth star:]  Illustrated in see Figure \ref{fig:geom} (a), the smooth star geometry is given by the parameterization 
 $\left(x(t),y(t)\right) = \left(1+0.3\cos(5t))\cos(t), \sin(t)(1+0.3\cos(5*t)\right)$ where $-\pi<t<\pi$.  
 \item[tear:] Illustrated in see Figure \ref{fig:geom} (b), the tear geometry is given by the parameterization 
 $\left(x(t),y(t)\right) = \left(2{\rm sign}(t)\sin(t),-\tan(\pi/4)\sin(2t)\right)$ where $-\pi<t<\pi$.  
 \item[pacman:] Illustrated in see Figure \ref{fig:geom} (c), the pacman geometry is given by the parameterization 
 $\left(x(t),y(t)\right) = \left({\rm sign}(t)\sin(1.5t),\tan(3\pi/2)\sin(t)\right)$ where $-\pi<t<\pi$.  
\end{list}

The Dirichlet-Neumann junctions occur at $t =0$ and $t=-\pi$ in 
 parameter space so that half the boundary has Dirichlet boundary conditions while the other half has 
 Neumann boundary conditions.  For each wavenumber/geometry combination, two types of boundary data are considered:
 
 \begin{list}{}{}
  \item[\textbf{known solution}:] Both the Dirichlet and Neumann boundary are generated by a known solution to the 
  Helmholtz problem corresponding to a collection of ten point charges inside of the geometry.  
  \item[\textbf{unknown solution}:] The Dirichlet boundary data is generated by the same ten point charges but the 
  Neumann data is set to zero.  Hence, the solution is not known a priori.
 \end{list}

For each wavenumber/geometry combination, tables \ref{tab:star_helm}-\ref{tab:pacman_helm} report 
\begin{list}{}{}
 \item[$l$:] The number of levels of refinement into the Dirichlet-Neumann junctions.
 \item[N:] The total number of discretization points.
 \item[$E_{\rm rel}$:] The relative error $E_{\rm rel} = \frac{\|u_l - u_{\rm ex}\|}{\|u_{\rm ex}\|}$ for the 
 known solution problem where 
 $u_l$ is the approximate solution, with $l$ levels of refinement, at twenty locations outside the geometry and $u_{\rm ex}$ is the 
 exact solution at the twenty locations outside the geometry generated by the ten interior point charges.
 \item[$E_{\rm conv}$:] The relative convergence error $E_{\rm conv} = \frac{\|u_l-u_{l-1}\|}{\|u_l\|}$ 
 where $u_l$ is the approximate solution, with $l$ levels of refinement, at twenty locations outside the geometry.
 \item[$\kappa(\mtx{A})$:] The condition number $\kappa(\mtx{A})$ of linear system resulting from the discretization 
 of (\ref{eq:sys}). 
\end{list}

% star geometry
\begin{table}[ht]
\centering
\begin{tabular}{|c|c|c|c|c|c|c|c|c|} \hline
 & \multicolumn{4}{|c|}{$\omega = 1$} & \multicolumn{4}{|c|}{$\omega = 10$}\\ \hline
$l$& $N$ & $E_{\rm rel}$ & $E_{\rm conv}$ & $\kappa(\mtx{A})$ & N & $E_{\rm rel}$ & $E_{\rm conv}$ & $\kappa(\mtx{A})$ \\ \hline
0 & 416&  7.05e-10& 1.00e-04& 3.06e03& 576 & 3.91e-10& 1.21e-04& 4.35e03 \\ \hline
1 & 480&  7.05e-10& 7.55e-05& 8.35e03& 640 & 3.91e-10& 9.04e-04& 1.13e04 \\ \hline
2 & 608&  7.05e-10& 2.36e-05& 4.30e04& 768 & 3.91e-10& 2.82e-05& 5.18e04 \\ \hline
3 & 864&  7.05e-10& 1.57e-06& 7.05e05& 1024& 3.91e-10& 1.88e-06& 8.38e05 \\ \hline
4 & 1376& 7.05e-10& 6.13e-09& 1.80e08& 1536& 3.91e-10& 7.36e-09& 2.14e08 \\ \hline
5 & 2400& 7.05e-10&   -     & 1.18e13& 2560& 3.91e-10&  -      & 1.40e13 \\ \hline
\end{tabular}\\
\begin{tabular}{|c|c|c|c|c|}\hline
 & \multicolumn{4}{|c|}{$\omega = 100$} \\ \hline
 $l$ & $N$&  $E_{\rm rel}$ & $E_{\rm conv}$ & $\kappa(\mtx{A})$ \\ \hline
 0 & 2656& 3.93e-10& 9.68e-05& 4.95e04\\ \hline
 1 & 2720& 3.93e-10& 7.25e-05& 1.01e05\\ \hline
 2 & 2848& 3.93e-10& 2.27e-05& 4.06e05\\ \hline
 3 & 3104& 3.93e-10& 1.51e-06& 6.49e06\\ \hline
 4 & 3616& 3.93e-10& 5.91e-09& 1.66e09\\ \hline
 5 & 4640& 3.93e-10& -       & 1.08e14\\ \hline
\end{tabular}
\caption{\label{tab:star_helm}The number of the levels of refinement $l$ into Dirichlet-Neumann junctions,
the number of discretization points $N$, the relative error $E_{\rm rel}$ for the boundary value problem
with known solution, the relative convergence error $E_{\rm conv}$ for the boundary value problem with 
unknown solution and the condition number $\kappa(\mtx{A})$ of the linear system resulting from the 
discretization of the integral equation (\ref{eq:sys}) on the smooth star geometry for three wave numbers $\omega = 1,10,$ and $100$. 
 }
 \end{table}

% tear geometry
\begin{table}[ht]
\centering
\begin{tabular}{|c|c|c|c|c|c|c|c|c|} \hline
 & \multicolumn{4}{|c|}{$\omega = 1$} & \multicolumn{4}{|c|}{$\omega = 10$}\\ \hline
$l$& $N$ & $E_{\rm rel}$ & $E_{\rm conv}$ & $\kappa(\mtx{A})$ & N & $E_{\rm rel}$ & $E_{\rm conv}$ & $\kappa(\mtx{A})$ \\ \hline
1 & 320&  4.48e-06& 2.93e-04& 3.17e03& 320 & 1.89e-07& 1.59e-03& 2.19e03 \\ \hline
2 & 448&  1.89e-06& 1.79e-04& 1.26e04& 488 & 4.37e-08& 8.82e-04& 8.70e03 \\ \hline
3 &704&  8.73e-08& 3.91e-05& 2.01e05& 704 & 2.58e-09& 1.62e-04& 1.39e05 \\ \hline
4 & 1216& 5.22e-10& 1.34e-06& 5.14e07& 1216& 1.19e-11& 4.52e-06& 3.56e07 \\ \hline
5 & 2240& 4.13e-11& 3.04e-08& 3.37e12& 2240& 7.37e-12& 7.027e-08& 2.34e12 \\ \hline
% 6 &     &         &         &        & 4288& 1.07e-07&   -      & 4.06e20 \\ \hline
\end{tabular}\\
\begin{tabular}{|c|c|c|c|c|}\hline
 & \multicolumn{4}{|c|}{$\omega = 100$} \\ \hline
 $l$ & $N$&  $E_{\rm rel}$ & $E_{\rm conv}$ & $\kappa(\mtx{A})$ \\ \hline
 1 & 2080& 1.43e-09& 1.11e-03& 3.35e04\\ \hline
 2 & 2208& 8.26e-10& 6.21e-04& 1.34e05\\ \hline
 3 & 2464& 7.67e-10& 1.13e-04& 2.15e06\\ \hline
 4 & 2976& 7.67e-10& 2.89e-06& 5.49e08\\ \hline
 5 & 4000& 7.67e-10& 2.74e-07& 3.55e13\\ \hline
%  6 & 6048& 1.88e-07&    -    &   \\ \hline
 \end{tabular}
\caption{\label{tab:tear_helm}The number of the levels of refinement $l$ into Dirichlet-Neumann junctions,
the number of discretization points $N$, the relative error $E_{\rm rel}$ for the boundary value problem
with known solution, the relative convergence error $E_{\rm conv}$ for the boundary value problem with 
unknown solution and the condition number $\kappa(\mtx{A})$ of the linear system resulting from the 
discretization of the integral equation (\ref{eq:sys}) on the tear geometry for three wave numbers $\omega = 1,10,$ and $100$. 
 }
 \end{table}

 % pacman geometry
\begin{table}[ht]
\centering
\begin{tabular}{|c|c|c|c|c|c|c|c|c|} \hline
 & \multicolumn{4}{|c|}{$\omega = 1$} & \multicolumn{4}{|c|}{$\omega = 10$}\\ \hline
$l$& $N$ & $E_{\rm rel}$ & $E_{\rm conv}$ & $\kappa(\mtx{A})$ & N & $E_{\rm rel}$ & $E_{\rm conv}$ & $\kappa(\mtx{A})$ \\ \hline
1 & 384&  2.23e-04& 2.93e-04& 4.93e03& 384 & 4.59e-04& 3.61e-04& 4.13e03 \\ \hline
2 & 512&  9.65e-05& 1.79e-04& 2.02e04& 512 & 1.95e-04& 1.79e-04& 1.69e04 \\ \hline
3 & 768&  1.74e-05& 3.91e-05& 3.23e05& 768 & 3.55e-05& 3.42e-05& 2.71e05 \\ \hline
4 & 1280& 5.82e-07& 1.34e-06& 8.28e07& 1280& 1.17e-06& 1.14e-06& 6.94e07 \\ \hline
5 & 2304& 1.39e-08& 3.04e-08& 5.43e12& 2304& 1.34e-09& 1.11e-07& 4.56e12 \\ \hline
% 6 & 4352& 8.61e-09&   -     & 2.11e20& 4352& 1.35e-07&     -   & 2.44e20 \\ \hline
\end{tabular}\\
\begin{tabular}{|c|c|c|c|c|}\hline
 & \multicolumn{4}{|c|}{$\omega = 100$} \\ \hline
 $l$ & $N$&  $E_{\rm rel}$ & $E_{\rm conv}$ & $\kappa(\mtx{A})$ \\ \hline
 1 & 2080& 5.71e-06& 8.59e-05& 5.91e04\\ \hline
 2 & 2208& 2.40e-06& 2.67e-05& 2.36e05\\ \hline
 3 & 2464& 4.38e-07& 1.79e-06& 3.79e06\\ \hline
 4 & 2976& 1.62e-08& 1.89e-08& 9.69e08\\ \hline
 5 & 4000& 7.06e-09& 1.96e-08& 6.26e13\\ \hline
%  6 & 6048& 1.76e-08& -       & -\\ \hline
\end{tabular}
\caption{\label{tab:pacman_helm}The number of the levels of refinement $l$ into Dirichlet-Neumann junctions,
the number of discretization points $N$, the relative error $E_{\rm rel}$ for the boundary value problem
with known solution, the relative convergence error $E_{\rm conv}$ for the boundary value problem with 
unknown solution and the condition number $\kappa(\mtx{A})$ of the linear system resulting from the 
discretization of the integral equation (\ref{eq:sys}) on the pacman geometry for three wave numbers $\omega = 1,10,$ and $100$. 
 }
 \end{table}

For the smooth star geometry (table \ref{tab:star_helm}), the proposed method is able to capture the smooth 
known solution without refinement independent of wavenumber.  The tear and pacman geometries have corners that
require refinement into the corners in order to capture the known solution since the boundary charge distribution $\sigma$
is not smooth in the corner.  Since $\sigma$ is smooth at the Dirichlet-Neumann junction away from the corner, 
the same performance would be observed if the refinement was only in the corner.  

For the unknown solution experiments, the boundary charge distribution is not smooth for any 
of the geometries.  Thus refinement into the Dirichlet-Neumann junctions is required to get 
high accuracy for all the experiments.  The performance of the proposed method is similar 
to its performance on the experiments with known solution on a geometry with a corner.  

The condition number of the linear system grows with the increased refinement into the 
Dirichlet-Neumann junctions.  This behavior is expected given the first kind nature of the 
integral equation.

\subsection{Laplace boundary value problem}
\label{sec:Lap}
 
This section reports on the performance of the proposed method and the integral formulation 
proposed in \cite{2009_Helsing_mixedBVP} for Laplace problems on the smooth star geometry 
with boundary data as specified by the {\it{known solution}} and {\it{unknown solution}} 
problems in the previous section.  It should be noted that \cite{2009_Helsing_mixedBVP} uses
a special quadrature to handle the hypersingular integral operator.  The results reported in 
this section use the same Nystr\"om discretization with 16-point composite Gaussian quadrature \cite{2001_rokhlin_kolm}
for both integral equations.  
The results using the integral formulation from \cite{2009_Helsing_mixedBVP} have an 
$H$ superscript.  

As with the Helmholtz boundary value problem, the proposed method does not require 
refinement when the boundary charge distribution is smooth.  When the boundary 
charge distribution is not smooth, refinement allows for technique to capture the 
solution to high accuracy.  This is contrast to the integral formulation from 
\cite{2009_Helsing_mixedBVP} which requires refinement to achieve high accuracy 
independent on whether or not the solution is smooth.  Note that while
the integral formulation from \cite{2009_Helsing_mixedBVP} is block second kind, since
the hypersingular term has not been dealt with specially, the condition number for 
the discretized linear system is nearly squared that of the first kind system.

\begin{table}[ht]
\centering
\begin{tabular}{|c|c|c|c|c|c|c|c|} \hline
l& N & $E_{\rm rel}$ & $\kappa(\mtx{A})$&$E^H_{\rm rel}$ & $\kappa(\mtx{A}^H)$ & $E_{\rm conv}$ & $E^H_{\rm conv}$ \\ \hline
0 & 352&  3.33e-10&9.21e02& 8.07e-05   & 6.78e04&  1.49e-06& 4.18e-03 \\ \hline
1 & 416&  3.33e-10&2.81e03& 4.04e-05   & 2.71e05&  1.12e-06& 1.62e-04 \\ \hline
2 & 544&  3.33e-10&1.28e04& 1.01e-05   & 4.33e06& 3.49e-07& 5.07e-05 \\ \hline
3 & 800& 3.33e-10& 2.12e05&6.31e-07    & 1.11e09&  2.32e-08& 3.37e-06 \\ \hline
4 & 1312& 3.33e-10& 5.44e07& 2.45e-09  &7.27e13 &  1.16e-11& 1.33e-08 \\ \hline
5 & 2336 & 3.33e-10 & 3.57e12& 1.69e-11& 3.12e23 &  - &   -       \\ \hline
\end{tabular}
\caption{\label{tab:Lap}The number of the levels of refinement $l$ into Dirichlet-Neumann junctions,
the number of discretization points $N$, the relative error $E_{\rm rel}$ for the Laplace boundary value problem
with known solution, the relative convergence error $E_{\rm conv}$ for the Laplace boundary value problem with 
unknown solution and the condition number $\kappa(\mtx{A})$ of the linear system resulting from the 
discretization of the integral equation (\ref{eq:sys}) on the star geometry. The values with the superscript $H$ 
correspond to the results from using the integral formulation from \cite{2009_Helsing_mixedBVP}.  }
 \end{table}

\section{Concluding remarks}

This paper presented a robust integral equation formulation for solving mixed boundary 
value problems of the form of (\ref{eq:mixedBVP}).  The formulation is a direct 
extension of the integral equation techniques for single boundary condition 
scattering problems.  Numerical results show 
that high accuracy can be obtained by utilizing local refinement near 
boundary condition junctions even if the junction is at a corner.   

If one does not have access to discretization techniques for hypersingular
kernels, the proposed solution technique is a high accuracy option for Laplace 
problems with mixed boundary value problems.

%%%%%%%%%%%%%%%%%%%%%%%%%%%%%%%%%%%%%%%%%%%%%%%%%%%%%%%%%%%%%%%%%%%%%%%%%%%%%%
\bibliography{main_bib}
\bibliographystyle{abbrv}

\end{document}